\documentclass[10pt,fleqn]{amsart}

\usepackage{amsmath,amssymb,latexsym}
\usepackage[mathscr]{eucal}
\usepackage{stmaryrd}

\theoremstyle{plain}
\newtheorem{theorem}{Theorem}
\newtheorem{lemma}{Lemma}

\newcounter{rmkctr}

\numberwithin{equation}{section}

\def\mydate{\number\year-\ifnum\month<10{0}\fi\number\month-\ifnum\day<10{0}\fi\number\day}

\newcommand{\dy}{\partial}

\newcommand{\ddt}[1]{\frac{\mathrm{d}{#1}}{\mathrm{d}{t}}}
\newcommand{\ddtau}[1]{\frac{\mathrm{d}{#1}}{\mathrm{d}{\tau}}}

\newcommand{\sfrac}[2]{{\textstyle\frac{#1}{#2}}}
\newcommand{\tssum}{{\textstyle\sum}}

\newcommand{\Zahl}{\mathbb{Z}}

\newcommand{\ex}{\mathrm{e}}
\newcommand{\im}{\mathrm{i}}
\newcommand{\eps}{\varepsilon}
\newcommand{\vfi}{\varphi}

\newcommand{\aand}{\quad\textrm{and}\quad}

\newcommand{\dtau}{\>\mathrm{d}\tau}

\newcommand{\ilapl}{\Delta^{-1}}

\newcommand{\sump}[1]{\mathop{\smash{\mathop{{\sum}_{#1}'}}{\vphantom\sum}}}

\newcommand{\xb}{{\boldsymbol{x}}}
\newcommand{\vb}{{\boldsymbol{v}}}

\newcommand{\fv}{f_{\vb}^{}}

\newcommand{\ZL}{\Zahl_L}
\newcommand{\jb}{{\boldsymbol{j}}}
\newcommand{\kb}{{\boldsymbol{k}}}
\newcommand{\lb}{{\boldsymbol{l}}}

\newcommand{\gb}{\nabla}
\newcommand{\sgb}{\nabla^\perp}

\newcommand{\divv}{\gb\!\cdot\!\vb}

\newcommand{\cnst}[1]{c_{#1}^{}}

\newcommand{\FD}{\mathsf{D}}

\newcommand{\BO}{B_\Omega}

\newcommand{\dstau}{\dy_\tau^*}

\newcommand{\Dom}{\mathscr{M}}
\newcommand{\Oh}{\mathsf{O}}

\newcommand{\cpoi}{c_0}

\newcommand{\wb}{\bar\omega}
\newcommand{\wt}{\tilde\omega}
\newcommand{\phb}{\bar\phi}
\newcommand{\pht}{\tilde\phi}
\newcommand{\fb}{\bar f}
\newcommand{\ft}{\tilde f}

\newcommand{\Lw}{\mathcal{L}}

\newcommand{\IO}{\mathsf{I}_\Omega}

\newcommand{\Attr}{\mathcal{A}}
\newcommand{\dima}{\textrm{dim}_\textrm{H}^{}\,\Attr}

\newcommand{\TB}{\mathcal{T}}

\newcommand{\supt}[1]{\llbracket{#1}\rrbracket}

\begin{document}

\title[Navier--Stokes on the $\beta$-plane]%
{Navier--Stokes equations on the $\beta$-plane}

\author[Al-Jaboori]{M.A.H.~Al-Jaboori}
\email{mustafa.hussain@durham.ac.uk}
\author[Wirosoetisno]{D.~Wirosoetisno}
\email{djoko.wirosoetisno@durham.ac.uk}
\urladdr{http://www.maths.dur.ac.uk/\~{}dma0dw}
\address{Department of Mathematical Sciences\\
   University of Durham\\
   Durham\ \ DH1~3LE, United Kingdom}

\thanks{This research was partially supported by ...}



\begin{abstract}
We show that, given a sufficiently regular forcing, the solution of
the two-dimensional Navier--Stokes equations on the periodic $\beta$-plane
(i.e.\ with the Coriolis force varying as $f_0+\beta y$) will become nearly
zonal: with the vorticity $\omega(x,y,t)=\wb(y,t)+\wt(x,y,t)$,
one has $|\wt|_{H^s}^2\le\beta^{-1} M_s(\cdots)$ as $t\to\infty$.
We use this show that, for sufficiently large $\beta$,
the global attractor of this system reduces to a point.
\end{abstract}

\maketitle


\section{Introduction}\label{s:intro}

The two-dimensional Navier--Stokes equations (2d NSE) have been the
subject of many studies and its basic mathematical properties (existence,
uniqueness, regularity, etc.) are now well understood;
see, e.g., \cite{doering-gibbon:aanse,temam:nsenfa} for reviews.
As a tool to understand various geophysical flows, it is often desirable
to include the effect of planetary rotation, but a constant rotation rate
(the so-called $f$-plane approximation) has no effect on the dynamics
when periodic boundary conditions are used.
To feel the effect of rotation, we need to go to the so-called
$\beta$-plane approximation, in which the rotation is given by
$f_0+\beta y$.

Simple physical arguments and numerical studies
\cite{rhines:75,maltrud-vallis:91} suggest that
a rotation rate that varies as $\beta y$ tends to force
the solution to become more zonal (a zonal flow is one that does not
depend on $x$).
In this article, we prove that this is indeed the case, by obtaining
a bound $|\wt(t)|_{H^s}^2\le \beta^{-1} M_s(f;\cdots)$, valid for large time $t$,
on the non-zonal part $\wt$ of the flow in terms of the forcing $f$.

With the further assumption that the forcing is independent of time,
it has been shown that the Navier--Stokes equations possess a global
attractor $\Attr$ of finite Hausdorff dimension.
The long-known and nearly optimal bound on this dimension
\cite{constantin-foias-temam:88} also applies to our rotating case,
but it does not take into account the effect of the rotation.
Using our bounds on $\wt$, we show that the dimension of $\Attr$
is zero for $\eps$ sufficiently small, reducing the long-time
dynamics to a single steady (and stable) flow determined completely
by the forcing.
This is to be contrasted with the situation for larger (but still small)
$\eps$, where the solution, although nearly zonal, evolves in time even
though $\dy_tf=0$.

Among the works similar in spirit to the present article, we mention
\cite{gallagher-straymond:07} where weak convergence to zonal flow is
proved for the (more difficult) $\beta$-plane shallow-water equations.
A related result for the inviscid Euler equation can be found in
\cite{tgs:87}.
The technique of using rapid oscillations to obtain better bounds have
been used in different contexts in, e.g., \cite{babin-al:aa:97,schochet:94}.

In the rest of this section, we describe the problem and set up the notation.
In Section~\ref{s:pre}, we review basic results on 2d NSE which will be
needed later.
The heart of this article is Section~\ref{s:bd0}, where $L^2$ and $H^s$
bounds are obtained for the non-zonal component of the flow.
An application of these bounds to the dimension of the global attractor
follows in Section~\ref{s:attr}.
The proof of an $L^\infty$ Agmon inequality is presented in the Appendix.

\medskip
In dimensional form, the two-dimensional Navier--Stokes equations read
\begin{equation}\label{q:fbeta}
   \dy_t\vb + \vb\cdot\gb\vb + \beta y\vb^\perp + \gb p
	= \mu\Delta\vb + \fv
\end{equation}
where the constant rotation $f_0$ has been dropped since it has
no effect (i.e.\ in 2d NSE, there is no difference between equatorial
and mid-latitude $\beta$-planes).
Here $\vb=(u,v)$ is the velocity with $\vb^\perp:=(-v,u)$
and $p$ is the pressure obtained by enforcing the incompressibility
constraint $\divv=0$.
In what follows, we will work with the dimensionless form
\begin{equation}\label{q:uv}\begin{aligned}
   &\dy_t\vb + \vb\cdot\gb\vb + \frac{Y}\eps \vb^\perp + \gb p
	= \mu\Delta\vb + \fv,\\
   &\divv = 0.
\end{aligned}\end{equation}

We work with $\xb=(x,y)\in\Dom:=[0,L_1]\times[-L_2/2,L_2/2]$, with periodic
boundary conditions in both directions.
Note that we have replaced $\beta y$ in \eqref{q:fbeta} by $Y(y)/\eps$,
where $Y(-L_2/2)=L_2/2$ and $Y(y)=y$ for $y\in(-L_2/2,L_2/2]$.
Furthermore, we assume the following symmetry on the velocity
\begin{equation}\label{q:sym}
   u(x,-y,t) = u(x,y,t)
   \aand
   v(x,-y,t) = -v(x,y,t).
\end{equation}
It is readily verified that if the initial data $\vb(0)$ and the forcing
$\fv(t)$ also satisfy this symmetry, which we henceforth assume,
it persists for all $t\ge0$.
Note also that periodicity and \eqref{q:sym} imply that
\begin{equation}\label{q:BCv}
   v(x,-L_2/2,t) = v(x,L_2/2,t) = 0.
\end{equation}
With no loss of generality, we require that the integral over $\Dom$
of $\vb$ vanishes.

In two dimensions, it is convenient to work with the vorticity
$\omega:=\sgb\!\cdot\!\vb=\dy_xv-\dy_yu$, whose evolution equation is
\begin{equation}\label{q:t00}
   \dy_t\omega + \vb\cdot\gb\omega + \frac{Y'}\eps\, v
	= \mu\Delta\omega + f.
\end{equation}
Here $f:=\sgb\!\cdot\!\fv$ and we can recover the velocity
using $\vb=\sgb\ilapl\omega$.
By our assumption on $\vb$, the integral of $\omega$ over $\Dom$ is zero;
similarly, $\ilapl$ is defined uniquely by the zero-integral condition.
The symmetry \eqref{q:sym} implies that $\omega(x,-y,t)=-\omega(x,y,t)$ and
\begin{equation}\label{q:BCw}
    \omega(x,-L_2/2,t) = \omega(x,L_2/2,t) = 0.
\end{equation}
Now $Y'(y)=1-L_2\,\delta(y-L_2/2)$, where $\delta$ is the Dirac distribution.
Using the fact that $v(x,\pm L_2/2,t)=0$, we replace $vY'$ by $v$
in \eqref{q:t00} and write
\begin{equation}\label{q:dwdt}
   \dy_t\omega + \vb\cdot\gb\omega + \frac1\eps\, v
	= \mu\Delta\omega + f.
\end{equation}
This is the form that we will be mostly working with.

It is also convenient to write \eqref{q:dwdt} in the usual functional form
\begin{equation}\label{q:abl}
   \dy_t\omega + B(\omega,\omega) + \frac1\eps L\omega + \mu A\omega = f,
\end{equation}
where $B(\omega,\omega^\sharp):=(\sgb\ilapl\omega)\cdot\gb\omega^\sharp$,
$L\omega:=\dy_x\ilapl\omega$ and $A:=-\Delta$.
The following properties, valid whenever the expressions make sense,
are readily verified by integration by parts and the boundary conditions
\eqref{q:BCv}--\eqref{q:BCw}
\begin{equation}\label{q:bla}\begin{aligned}
   &(B(\omega,\omega^\sharp),\omega^\sharp)_{L^2}^{} = 0,\\
   &(L\omega,\omega)_{L^2}^{} = 0,\\
   &(A\omega,\omega)_{L^2}^{} = |\gb\omega|_{L^2}^2.
\end{aligned}\end{equation}


\section{Preliminary Estimates}\label{s:pre}

The estimates derived in this section are standard from the theory of 2d NSE
(see, e.g., \cite{temam:iddsmp,doering-gibbon:aanse,robinson:idds}),
with very minor modifications to handle the Coriolis term.
We gather them here for later use.

We start by noting that the vanishing of spatial integrals of $\vb$
and $\omega$ implies the equivalence of the norms $|\omega|_{H^s}^{}$ and
$|\gb^s\omega|_{L^2}^{}:=|(-\Delta)^{s/2}\omega|_{L^2}^{}$,
which will thus be used interchangeably below.
We denote by $\cpoi$ the constant in Poincar{\'e} inequality
\begin{equation}\label{q:poidef}
   \cpoi\,|\gb^s\omega|_{L^2}^{} \le |\gb^{s+1}\omega|_{L^2}^{}.
\end{equation}
Besides the usual Sobolev and interpolation inequalities for two and
one dimensions (for functions depending on $y$ only), we note the
one-dimensional Agmon inequality
\begin{equation}
   |\bar w|_{L^\infty}^{} \le c\,|\bar w|_{L^2}^{1/2}|\gb\bar w|_{L^2}^{1/2}.
\end{equation}
A version we use for the two-dimensional case is in Appendix~\ref{s:linfty}.

The $L^2$ estimate for the velocity is obtained by multiplying \eqref{q:uv}
by $\vb$ and using Cauchy--Schwarz,
\begin{equation}\label{q:bdv}
   \ddt{}|\vb|_{L^2}^2 + \mu\,|\gb\vb|_{L^2}^2 \le \frac{c}\mu\,|\fv|_{L^2}^2\,.
\end{equation}
Assuming that $\fv\in L_t^\infty L_{\xb}^2$, we thus have
$\vb\in L_t^\infty L_{\xb}^2 \cap L_t^2 H_{\xb}^1$.
Here and henceforth, $L_t^pL_{\xb\vphantom{t}}^q:=L^p((0,\infty);L^q(\Dom))$, and
$H_{\xb}^s$ and $|\gb^s\omega|$ below are defined in the usual way.
We denote
\begin{equation}\label{q:kmdef}
   \supt{w} := \sup_{t>0}\,|w(t)|_{L^2}^{}.
\end{equation}
Here and elsewhere in this article, $c$ denotes a generic constant depending
only on $\Dom$ whose value may not be the same each time it appears,
while numbered constants such as $c_1^{}$ have fixed values.

Now let $\vfi\in L^\infty$ be such that $\vfi'(t)=1$ for $t\in[0,\sfrac12]$
and $\vfi(t)=1$ for $t>2$, so $\vfi(t)\simeq\tanh t$.
Multiplying \eqref{q:dwdt} by $\vfi\omega$ in $L^2$, or equivalently,
multiplying (\ref{q:uv}a) by $-\vfi\Delta\vb$ in $L^2$, we obtain
\begin{equation}\label{q:bd0}
   \ddt{}\bigl(\vfi|\omega|_{L^2}^2\bigr) + \mu\,\vfi\,|\gb\omega|_{L^2}^2
	\le \vfi'\,|\omega|^2 + \frac{c}\mu\,\vfi\,|\fv|_{L^2}^2\,.
\end{equation}
Assuming henceforth that
$\supt{\fv}<\infty$, we have
$\omega\in L_t^\infty L_{\xb}^2 \cap L_t^2 H_{\xb}^1$ and,
for $t\ge T_0(|\vb(0)|_{L^2}^{},\supt{\fv};\mu)$,
\begin{equation}\label{q:bd0u}
   |\omega(t)|_{L^2}^{} \le \frac{c}{\mu}\,\supt{\fv}.
\end{equation}
Note that $T_0$ does not depend on $\omega(0)$ and that
the requirement $\fv\in L_t^\infty L_{\xb}^2$ can be
weakened in $t$, but we shall not do so here.

A bound in $H^m$ is obtained as follows.
Fix a multi-index $\alpha=(\alpha_1,\alpha_2)$ with
$|\alpha|:=\alpha_1+\alpha_2=m$, and multiply \eqref{q:dwdt} by
$\FD^{2\alpha}\omega:=\dy_x^{2\alpha_1}\dy_y^{2\alpha_2}\omega$ in $L^2$,
\begin{equation}
   \frac12\ddt{\;}|\FD^\alpha\omega|^2 + (\vb\cdot\gb\omega,\FD^{2\alpha}\omega) + \frac1\eps (v,\FD^{2\alpha}\omega) + \mu\,|\gb\FD^\alpha\omega|^2 = (f,\FD^{2\alpha}\omega)
\end{equation}
where here and henceforth $|\cdot|$ and $(\cdot,\cdot)$ denote
$L^2$ norm and inner product.
The linear term involving $1/\eps$ vanishes, and one then proceeds as usual:
Using the fact that $(\vb\cdot\gb\FD^\alpha\omega,\FD^\alpha\omega)=0$, the
nonlinear term is bounded as
\begin{equation}\label{q:whs}\begin{aligned}
   \bigl|(\vb\cdot\gb\omega,\FD^{2\alpha}\omega)\bigr|
	&\le \sum_{1\le|\beta|\le|\alpha|} \bigl|((\FD^\beta\vb)\cdot\gb \FD^{\alpha-\beta}\omega,\FD^\alpha\omega)\bigr|\\
	&\le c\,\tssum_\beta\,\bigl|\FD^{\beta-1}\omega\bigr|_{L^4}^{} \bigl|\FD^{\alpha-\beta+1}\omega\bigr|_{L^4}^{} \bigl|\FD^\alpha\omega\bigr|_{L^2}^{}\\
	&\le c\,\tssum_\beta\,\bigl|\FD^{\beta-1}\omega\bigr|_{H^{1/2}}^{} \bigl|\FD^{\alpha-\beta+1}\omega\bigr|_{H^{1/2}}^{} \bigl|\FD^\alpha\omega\bigr|_{L^2}^{}\\
	&\le c(m)\,\tssum_{l=1}^m\,|\omega|_{H^{l-1/2}}^{}|\omega|_{H^{m-l+3/2}}^{}|\omega|_{H^m}^{}
\end{aligned}\end{equation}
where we have used Sobolev inequalities for the second and
third line, and where $l:=|\beta|$ in the last line.
Using the interpolation inequalities
\begin{equation}\label{q:interp0}\begin{aligned}
   &|\omega|_{H^{l-1/2}}^{} \le c\,|\omega|_{L^2}^{(2m-2l+3)/(2m+2)} |\omega|_{H^{m+1}}^{(2l-1)/(2m+2)}\\
   &|\omega|_{H^{m-l+3/2}}^{} \le c\,|\omega|_{L^2}^{(2l-1)/(2m+2)} |\omega|_{H^{m+1}}^{(2m-2l+3)/(2m+2)},
\end{aligned}\end{equation}
followed by Cauchy--Schwarz and summing over $\alpha$, we obtain
\begin{equation}\label{q:bdwm0}
   \ddt{\;}|\omega|_{H^m}^2 + \frac{3\mu}{2}\,|\omega|_{H^{m+1}}^2
	\le \frac{c(m)}\mu\,|\omega|_{L^2}^2|\omega|_{H^m}^2 + \frac{c'(m)}\mu\,|f|_{H^{m-1}}^2
\end{equation}
for $m=1,2,\cdots$.
Proceeding by Gronwall and induction on \eqref{q:bdwm0},
we have the following uniform bounds independent of the initial data
\begin{equation}\label{q:bdwm}\begin{aligned}
   &|\omega(t)|_{H^m}^2 \le \Bigl(\frac{c(m)}\mu\,\supt{\gb^{m-1}f}^2 + 1\Bigr)^{m+1}\\
   &\ex^{-\nu t'} \int_t^{t+t'} \ex^{\nu\tau} |\omega(\tau)|_{H^{m+1}}^2 \dtau \le \frac{c}\mu\,\Bigl(\frac{c(m)}\mu\,\supt{\gb^{m-1}f}^2 + 1\Bigr)^{m+1}
\end{aligned}\end{equation}
valid for all $t\ge T_m(|\vb(0)|_{L^2}^{},\supt{\gb^{m-1}f}^{};\mu)$.
Here $\nu=\cpoi^2\mu$ with $\cpoi^{}$ the constant in Poincar{\'e} inequality
\eqref{q:poidef}.
Note that for $T_m$ to depend only on $|\vb(0)|_{L^2}^{}$, and for the validity
of \eqref{q:bdwm} for $\vb(0)\not\in H^1$, we need to multiply
by $\vfi$ as in \eqref{q:bd0}, but this was not done explicitly for conciseness.


\section{Bounds on the Non-zonal Component}\label{s:bd0}

Assuming sufficient regularity for $f$, which implies that for $\omega$
for any $t>0$, we expand them in Fourier series
\begin{equation}\label{q:wkdef}\begin{aligned}
   &\omega(\xb,t) = \tssum_\kb\,\omega_\kb(t)\, \ex^{\im\kb\cdot\xb-\im\Omega_\kb t/\eps}\\
   &f(\xb,t) = \tssum_\kb\,f_\kb(t)\, \ex^{\im\kb\cdot\xb}
\end{aligned}\end{equation}
where $\kb=(k_1,k_2)\in\ZL:=\{(2\pi l_1/L_1,2\pi l_2/L_2): (l_1,l_2)\in\Zahl^2\}$
and $\Omega_\kb:=-k_1/|\kb|^2$ is i times the eigenvalue of
the linear operator $L$ for wavenumber $\kb$.
Since $\omega$ and $f$ have vanishing integrals over $\Dom$,
$\omega_{\boldsymbol{0}}=0$ and $f_{\boldsymbol{0}}=0$.
Here and in what follows, sums over wavenumbers are understood to be taken
over $\ZL$.
In terms of Fourier components, \eqref{q:dwdt} reads
\begin{equation}\label{q:dwkdt}
   \ddt{\omega_\lb} + \tssum_{\jb\kb}\,B_{\jb\kb\lb}\omega_\jb\omega_\kb\,\ex^{\im(\Omega_\lb-\Omega_\jb-\Omega_\kb)t/\eps}
	+ \mu\,|\lb|^2\omega_\lb = f_\lb\,\ex^{\im\Omega_\lb t/\eps}
\end{equation}
where the coefficient of the nonlinear term is
\begin{equation}\label{q:bjkldef}
   B_{\jb\kb\lb} = (B(\ex^{\im\jb\cdot\xb},\ex^{\im\kb\cdot\xb}),\ex^{\im\lb\cdot\xb})
	= |\Dom|\,\frac{\jb\wedge\kb}{|\jb|^2}\,\delta_{\jb+\kb-\lb}
\end{equation}
with $\jb\wedge\kb:=j_1k_2-j_2k_1$.
We note that the linear term $\eps^{-1}L\omega$ has been removed
from \eqref{q:dwkdt} by including $\exp(-\im\Omega_\kb t/\eps)$
in (\ref{q:wkdef}a).

Let us split $\omega$ into a slow part $\wb$, for which $\Omega_\kb=0$,
and the remaining fast part $\wt:=\omega-\wb$, viz.,
\begin{equation}\begin{aligned}
   &\wb(\xb,t) = \tssum_{k_1=0}\,\omega_\kb(t)\, \ex^{\im\kb\cdot\xb}\\
   &\wt(\xb,t) = \tssum_{k_1\ne0}\,\omega_\kb(t)\, \ex^{\im\kb\cdot\xb-\im\Omega_\kb t/\eps}.
\end{aligned}\end{equation}
We note that, also having zero integrals over $\Dom$,
$\wt$ and $\wb$ are orthogonal in $H^m$ for $m=0,1,\cdots$.
For convenience, we also define
\begin{equation}
   \wb_\kb := \begin{cases} \omega_\kb &\textrm{if } k_1=0\\ 0 &\textrm{otherwise,}\end{cases}
   \qquad\textrm{and}\qquad
   \wt_\kb := \begin{cases} \omega_\kb &\textrm{if } k_1\ne0\\ 0 &\textrm{otherwise.}\end{cases}
\end{equation}
Our objective in this section is to obtain long-time bounds for $\wt$
that tend to zero as $\eps\to0$.


\subsection{Bound in $L^2$}\label{s:l2}

The development in this subsection largely follows that in
\cite{temam-dw:lbal} for the primitive equations, the main difference
being the absence of a spectral gap (that is, the eigenvalues of the
antisymmetric operator $L$ accumulate at zero in the present case).

We start by multiplying \eqref{q:abl} by $\wt$ in $L^2$,
\begin{equation}\label{q:dwtdt0}
   (\dy_t\omega,\wt) + (B(\omega,\omega),\wt) + \frac1\eps (L\omega,\wt)
	+ \mu\,(A\omega,\wt) = (f,\wt).
\end{equation}
Now, using (\ref{q:bla}a) twice and the fact that $B(\wb,\wb)=0$,
\begin{equation}\begin{aligned}
   (B(\omega,\omega),\wt) &= (B(\omega,\wt),\wt) + (B(\omega,\wb),\wt)\\
	&= (B(\wb,\wb),\wt) + (B(\wt,\wb),\wt)\\
	&= - (B(\wt,\wt),\wb).
\end{aligned}\end{equation}
Thus \eqref{q:dwtdt0} becomes
\begin{equation}\label{q:dwtdt1}
   \frac12\ddt{\;}|\wt|^2 + \mu\,|\gb\wt|^2 = (B(\wt,\wt),\wb) + (f,\wt).
\end{equation}
Dropping the nonlinear term for the moment, the fact that $\wt$ is
rapidly varying while $f$ is slow implies that the effective forcing
from the rhs becomes weaker for smaller $\eps$.
This essentially is the mechanism for the attenuation of the fast part $\wt$;
the nonlinear term will be handled in the proof below.
Recalling the definition \eqref{q:kmdef},
we state the result of this subsection.

\begin{theorem}\label{t:o1}
Assume that the initial data $\vb(0)\in L^2(\Dom)$
and that the forcing is bounded as $\supt{\gb^2f}+\supt{\dy_tf}<\infty$.
Then there exist $\TB_0(|\vb(0)|_{L^2}^{},\supt{\gb^2f},\supt{\dy_tf};\mu)$ and
$M_0(\supt{\gb^2f},\supt{\dy_tf};\mu)$ such that, for $t\ge \TB_0$,
\begin{equation}\label{q:bdwt}\begin{aligned}
   &|\wt(t)|_{L^2}^2 \le \eps\,M_0,\\
   &\mu\,\ex^{-\nu (t+t')}\int_t^{t+t'} \ex^{\nu\tau} |\gb\wt(\tau)|^2 \dtau
	\le \eps\,M_0.
\end{aligned}\end{equation}
\end{theorem}


\noindent{\sc Proof.}
Recalling that $\nu=\mu\cpoi^2$, we obtain from \eqref{q:dwtdt1}
\begin{equation}\label{q:dwtdt}
   \ddt{\;}\bigl(\ex^{\nu t}|\wt|^2\bigr) + \mu\ex^{\nu t}|\gb\wt|^2
	\le 2\ex^{\nu t}(B(\wt,\wt),\wb) + 2\ex^{\nu t}(f,\wt).
\end{equation}

We integrate the last term from $0$ to $t$ by parts,
\begin{equation}\begin{aligned}
   \int_0^t \ex^{\nu\tau} (f,\wt)\dtau
	&= |\Dom|\, \tssum_\kb' \int_0^t f_\kb(\tau)\overline{\wt_\kb(\tau)}\,
		\ex^{\im\Omega_\kb\tau/\eps+\nu\tau} \dtau\\
	&= -\im\eps|\Dom|\,\sump{\kb}\, \frac1{\Omega_\kb}\Bigl[ f_\kb\overline{\wt_\kb}\,\ex^{\im\Omega_\kb\tau/\eps+\nu\tau}\Bigr]_0^t\\
	&\qquad {}+ \im\eps|\Dom|\, \sump{\kb}\, \frac1{\Omega_\kb} \int_0^t \ddtau{\;}\bigl[ f_\kb\overline{\wt_\kb}\,\ex^{\nu\tau}\bigr] \ex^{\im\Omega_\kb\tau/\eps} \dtau
\end{aligned}\end{equation}
where the prime on the sums indicates that the resonant terms (i.e.\ those
with $\Omega_\kb=0$) are excluded.
Defining the operator $\dy_t^*$ by, for any $w$ for which it makes sense,
\begin{equation}\label{q:dstaudef}\begin{aligned}
   &\dy_t^* w := \ex^{-tL/\eps}\dy_t\bigl(\ex^{tL/\eps}w\bigr)\\
   \Rightarrow\qquad
   &\dy_t^*\wt := \dy_t\wt + \frac1\eps L\wt
	= -\tilde B(\omega,\omega) - \mu A\wt + \ft,
\end{aligned}\end{equation}
and defining the operator $\IO$ by
\begin{equation}\label{q:IOdef}
   \IO\ft(\xb,t)
	:= \sump{\kb} \frac1{\im\Omega_\kb}\, f_\kb(t)\,\ex^{\im\kb\cdot\xb}
	= \im \sump{\kb} \frac{|\kb|^2}{k_1}\, f_\kb(t)\,\ex^{\im\kb\cdot\xb},
\end{equation}
which being the restricted inverse of $L$ is also antisymmetric, we can write
\begin{equation}\label{q:i00}\begin{aligned}
   \int_0^t \ex^{\nu\tau} (f,\wt)\dtau
	&= \eps\,(\IO\ft,\wt)(t)\ex^{\nu t} - \eps\,(\IO\ft,\wt)(0)\\
	&- \eps \int_0^t \bigl[ \nu(\IO\ft,\wt) + (\dy_\tau\IO\ft,\wt)
		+ (\IO\ft,\dstau\wt)\bigr] \ex^{\nu\tau} \dtau.
\end{aligned}\end{equation}
Using \eqref{q:IOdef}, the endpoint terms can be bounded as
\begin{equation}\label{q:j04}
   |(\IO\ft,\wt)| \le c\,|\gb\ft|\,|\gb\wt|.
\end{equation}
We now bound the terms in the integrand.
First,
\begin{equation}
   \bigl|(\dy_\tau\IO\ft,\wt)\bigr|
	= \bigl|(\dy_\tau\ft,\IO\wt)\bigr|
	\le c\,|\dy_\tau\ft|\,|\Delta\wt|.
\end{equation}
Next, using (\ref{q:dstaudef}b) and noting
the fact that $(\IO\ft,\ft)=0$, we bound the last term in \eqref{q:i00} by
\begin{equation}
   \bigl|(\IO\ft,\mu\Delta\wt)\bigr|
	\le \mu c\,|\Delta\ft|\,|\Delta\wt|;
\end{equation}
and, using Sobolev and interpolation inequalities,
\begin{equation}\begin{aligned}
   \bigl|(\IO\ft,B(\omega,\omega))\bigr|
	&\le c\,|\gb\ft|_{L^2}^{}|\gb B(\wt,\wt)|_{L^2}^{}\\
   &\le c\,|\gb\ft|_{L^2}^{} |\omega|_{L^4}^{}|\gb\omega|_{L^4}^{}
	+ c\,|\gb\ft|_{L^2}^{} |\gb^{-1}\omega|_{L^\infty}^{}|\Delta\omega|_{L^2}^{}\\
	&\le c\,|\gb\ft|\,|\gb\omega|\,|\Delta\omega|.
\end{aligned}\end{equation}
Thus the integral in \eqref{q:i00} is bounded as
\begin{equation}\label{q:j02}\begin{aligned}
   \biggl|\int_0^t &[\nu(\IO\ft,\wt) + (\dy_\tau\IO\ft,\wt)
		+ (\IO\ft,\dstau\wt)]\,\ex^{\nu\tau} \dtau\biggr|\\
	&\le c \int_0^t \bigl[ \mu\,|\Delta\ft|\,|\Delta\wt|
		+ |\dy_\tau\ft|\,|\Delta\wt| + |\gb\ft|\,|\gb\omega|\,|\Delta\omega| \bigr] \ex^{\nu\tau}\dtau\\
	&\le c\int_0^t \bigl\{ (1+\mu)\,|\Delta\wt|^2 + \mu\,|\Delta\ft|^2 + |\dy_\tau\ft|^2 + |\Delta\omega|\,|\gb\omega|\,|\gb\ft| \bigr\}\,\ex^{\nu\tau}\dtau.
\end{aligned}\end{equation}

We now treat the penultimate term in \eqref{q:dwtdt}.
First, we write
\begin{equation}\label{q:Bsym}\begin{aligned}
   (B(\wt,\wt),\wb) &= \tssum_{\jb\kb\lb}\,B_{\jb\kb\lb}\wt_\jb\wt_\kb\overline{\wb_\lb}\,\ex^{-\im(\Omega_\jb+\Omega_\kb)t/\eps}\\
	 &= \frac12 \tssum_{\jb\kb\lb}\,(B_{\jb\kb\lb}+B_{\kb\jb\lb})\,\wt_\jb\wt_\kb\overline{\wb_\lb}\,\ex^{-\im(\Omega_\jb+\Omega_\kb)t/\eps}
\end{aligned}\end{equation}
and then note that $B_{\jb\kb\lb}+B_{\kb\jb\lb}=0$ in the resonant case,
i.e.\ when $\Omega_\jb+\Omega_\kb=0$ and $l_1=0$.
Furthermore, we have
\begin{equation}\label{q:nores}\begin{aligned}
   B_{\jb\kb\lb} + B_{\kb\jb\lb}
	&= \biggl(\frac{\jb\wedge\kb}{|\jb|^2} + \frac{\kb\wedge\jb}{|\kb|^2}\biggr)\,|\Dom|
	& &= (\jb\wedge\kb)\,\biggl(\frac1{|\jb|^2} - \frac1{|\kb|^2}\biggr)\,|\Dom|\\
	&= j_1l_2\,\biggl(\frac1{|\jb|^2} - \frac1{|\kb|^2}\biggr)\,|\Dom|
	& &= -l_2\,(\Omega_\jb+\Omega_\kb)\,|\Dom|
\end{aligned}\end{equation}
whenever $\jb+\kb=\lb$ and $l_1=0$.
Motivated by \eqref{q:Bsym}, we introduce the bilinear symmetric
operator $\BO$ by
\begin{equation}\label{q:BOdef}\begin{aligned}
   (\BO(\wt^\sharp,\wt^\flat),\wb) &:=
	\frac{\im}{2} \sump{\jb\kb\lb}\, \frac{B_{\jb\kb\lb}+B_{\kb\jb\lb}}{\Omega_\jb+\Omega_\kb}\,\wt_\jb^\sharp\wt_\kb^\flat\overline{\wb_\lb}\,\ex^{-\im(\Omega_\jb+\Omega_\kb)t/\eps}\\
	&= \frac{|\Dom|}{2\im} \sump{\jb\kb\lb}\, l_2 \,\wt_\jb^\sharp\wt_\kb^\flat\overline{\wb_\lb}\,\ex^{-\im(\Omega_\jb+\Omega_\kb)t/\eps}
\end{aligned}\end{equation}
for any $\wt^\sharp$, $\wt^\flat$ and $\wb$, where the prime on the sum
again indicates that resonant terms (for which $\Omega_\jb+\Omega_\kb=0$)
are omitted.
We note that, thanks to \eqref{q:nores}, the resonant terms are also absent
in $(B(\wt,\wt),\wb)$.
Integrating by parts, we have
\begin{equation}\label{q:j01}\begin{aligned}
   \int_0^t \ex^{\nu\tau}(B(\wt,\wt),\wb) \dtau
	&= \eps\,\ex^{\nu t} (\BO(\wt,\wt),\wb)(t) - \eps\,(\BO(\wt,\wt),\wb)(0)\\
	&+ \eps \int_0^t \ex^{\nu\tau} \bigl[ \nu (\BO(\wt,\wt),\wb)
		+ 2\,(\BO(\dstau\wt,\wt),\wb)\\
		&\hbox to116pt{}+ (\BO(\wt,\wt),\dy_\tau\wb) \bigr] \dtau.
\end{aligned}\end{equation}
For the last term in the integrand, we use the fact that
$\bar B(\wt,\wb) = \bar B(\wb,\wt) = \bar B(\wb,\wb) = 0$ to write
\begin{equation}
   \dy_\tau\wb = -\bar B(\wt,\wt) - \mu A\wb + \fb
\end{equation}
and estimate, using $H^1\subset L^\infty$ for $\fb$ and \eqref{q:agmon}
for the $L^\infty$ estimates,
\begin{equation}\begin{aligned}
   \bigl|&(\BO(\wt,\wt),\dy_\tau\wb)\bigr|
	\le c\,|\wt|_{L^2}^{}|\dy_y\wt|_{L^2}^{}|\fb|_{L^\infty}^{}
	+ \mu c\,|\wt|_{L^4}^{}|\dy_y\wt|_{L^4}^{}|\Delta\wb|_{L^2}^{}\\
	&\hbox to160pt{} {}+ c\,|\wt|_{L^\infty}^{}|\dy_y\wt|_{L^2}^{}|\gb^{-1}\wt|_{L^\infty}^{}|\gb\wt|_{L^2}^{}\\
	&\le c\,|\wt|\,|\gb\wt|\,|\fb'|
	+ \mu c\,|\wt|^{1/2}|\gb\wt|\,|\Delta\omega|^{3/2}+ c\,|\gb\omega|^3|\wt|\,\Bigl(\log\frac{|\Delta\omega|}{\cpoi|\gb\omega|}+c'\Bigr).
\end{aligned}\end{equation}
For the term involving $(\BO(\dstau\wt,\wt),\wb)$, we bound,
using $\dstau\wt+\tilde B(\omega,\omega) + \mu A\wt = \ft$
and the inequality $|\wb|_{L^\infty}^{}\le c\,|\wb|^{1/2}|\wb'|^{1/2}$,
\begin{equation}\begin{aligned}
   \bigl|(\BO(\dstau\wt,\wt),\wb)\bigr|
	&\le c\,|\dy_y\ft|\,|\wt|\,|\wb|_{L^\infty}^{} + c\,|\ft|\,|\dy_y\wt|\,|\wb|_{L^\infty}^{}\\
	&\qquad {}+ \mu c\,|\Delta\omega|\,|\wt|\,|\wb'|_{L^\infty}^{}
	+ c\,|\gb^{-1}\omega|_{L^\infty}^{}|\gb\omega|\,|\wt|\,|\wb'|_{L^\infty}^{}\\
	&\le c\,|\gb\ft|\,|\wt|\,|\wb'| + c\,|\ft|\,|\gb\wt|\,|\wb'|
	+ \mu c\,|\wt|\,|\wb'|^{1/2}|\Delta\omega|^{3/2}\\
	&\qquad {}+ c\,|\omega|^2|\gb\omega|^{3/2}|\wb''|^{1/2}\Bigl(\log\frac{|\gb\omega|}{\cpoi|\omega|}+1\Bigr)^{1/2}
\end{aligned}\end{equation}
where all unadorned norms are $L^2$.
Finally, we bound
\begin{equation}\begin{aligned}
   \bigl|(\BO(\wt,\wt),\wb)\bigr|
	&\le c\,|\wt|\,|\dy_y\wt|\,|\wb|_{L^\infty}^{}\\
	&\le c\,|\wt|\,|\gb\wt|\,|\wb|^{1/2}|\wb'|^{1/2}.
\end{aligned}\end{equation}
Using these also to bound the endpoint terms,
the integral in \eqref{q:j01} is bounded as
\begin{equation}\label{q:j03}\begin{aligned}
   \biggl|\int_0^t \ex^{\nu\tau}&(B(\wt,\wt),\wb) \dtau\biggr|\\
	&\le \eps c\,\bigl[|\wt|\,|\gb\wt|\,|\wb|^{1/2}|\wb'|^{1/2}\bigr](t)\,\ex^{\nu t}
	+ \eps c\,\bigl[|\wt|\,|\gb\wt|\,|\wb|^{1/2}|\wb'|^{1/2}\bigr](0)\\
	& {}+ \eps \int_0^t \Bigl\{
		c\,|\gb f|\,|\wt|\,|\gb\omega| + c\,|\ft|\,|\gb\omega|^2
		+ \mu c\,|\omega|^{1/2}\,|\gb\omega|\,|\Delta\omega|^{3/2}\\
	&\hbox to75pt{} {}+ c\,|\omega|\,|\gb\omega|^{5/2}|\Delta\omega|^{1/2}\Bigl(\log\frac{|\Delta\omega|}{\cpoi|\gb\omega|}+c'\Bigr) \Bigl\}\,\ex^{\nu\tau} \dtau.
\end{aligned}\end{equation}

Putting together \eqref{q:j04}, \eqref{q:j02} and \eqref{q:j03}, we have
\begin{equation}\label{q:bd1z}\begin{aligned}
   |\wt(&t)|^2 + \mu \int_0^t |\gb\wt|^2 \ex^{\nu(\tau-t)} \dtau
	\le \ex^{-\nu t}|\wt(0)|^2\\
	&+ \eps \cnst2\,(1+\ex^{-\nu t})\,\sup_{0\le t'\le t}\,\bigl\{
		|\gb\ft|\,|\gb\wt| + |\omega|^{3/2}|\gb\omega|^{3/2} \bigr\}\\
	&+ \eps \cnst3(\mu) \int_0^t \Bigl\{ |\Delta\ft|^2 + |\dy_\tau\ft|^2
		+ |\gb\ft|\,|\gb\omega|\,|\Delta\omega| + |\Delta\omega|^2\bigl(1+|\gb\omega|\bigr)\\
	&\hbox to80pt{} {}+ |\omega|\,|\gb\omega|^{5/2}|\Delta\omega|^{1/2}\Bigl(\log\frac{|\Delta\omega|}{\cpoi|\gb\omega|}+c'\Bigr) \Bigl\}\,\ex^{\nu(\tau-t)} \dtau.
\end{aligned}\end{equation}
We now shift the origin of time such that $t=0$ corresponds to $T_2$
in \eqref{q:bdwm}.
The hypothesis that $\supt{\gb^2f}+\supt{\dy_tf}<\infty$ then implies that
both the endpoints and the integral in \eqref{q:bd1z} are bounded uniformly
for all $t>0$, independently of the initial data provided that
$\vb\in L^2$ initially.
Rewriting the bound in \eqref{q:bd1z} as
\begin{equation}\label{q:bdl2}
   |\wt(t)|^2 + \mu \int_0^t |\gb\wt|^2 \ex^{\nu(\tau-t)} \dtau
	\le \ex^{-\nu t}|\wt(0)|^2
	+ \frac{\eps}2\,M_0(\supt{\gb^2f},\supt{\dy_tf};\mu),
\end{equation}
the proof is complete.

We also note from \eqref{q:bd1z} that the hypothesis
$f\in L_t^\infty H_{\xb}^2$ and $\dy_t f\in L_t^\infty L_{\xb}^2$
can be weakened to
$f\in L_t^2 H_{\xb}^2\cap L_t^\infty H_{\xb}^1$ and $\dy_t f\in L_t^2 L_{\xb}^2$.




\subsection{Bounds in $H^s$}\label{s:hs}
With a little extra work, $H^s$ bounds for $\wt$ that scales as $\sqrt\eps$
can also be obtained.
We do this explicitly for $|\gb\wt|$ and sketch the computation
for $s=2,3,\cdots$.

\medskip
For the $H^1$ bound, we multiply \eqref{q:abl} by $A\wt$ in $L^2$ to get
\begin{equation}
   \frac12\ddt{\;}|\gb\wt|^2 + \mu\,|\Delta\wt|^2
	+ (B(\omega,\omega),A\wt) = (\ft,A\wt),
\end{equation}
which implies [cf.\ \eqref{q:dwtdt}]
\begin{equation}
   \ddt{\;}\bigl(\ex^{\nu t}|\gb\wt|^2\bigr) + \mu\,\ex^{\nu t}|\Delta\wt|^2
	\le 2\ex^{\nu t}\,(B(\omega,\omega),\Delta\wt)
	- 2\ex^{\nu t}(\ft,\Delta\wt).
\end{equation}
As in the $L^2$ case, we integrate from 0 to $t$,
\begin{equation}\begin{aligned}
   \ex^{\nu t}|\gb\wt(t)|^2 - |\gb\wt(0)|^2
	&+ \mu\! \int_0^t |\Delta\wt|^2 \,\ex^{\nu\tau}\dtau\\
	&\le 2 \int_0^t \bigl\{ (B(\omega,\omega),\Delta\wt) - (\ft,\Delta\wt) \bigr\}\, \ex^{\nu\tau} \dtau.
\end{aligned}\end{equation}
The forcing term gives
\begin{equation}\label{q:j00}\begin{aligned}
   \int_0^t \ex^{\nu\tau} (f,\Delta\wt)\dtau
	&= \eps\,(\IO\ft,\Delta\wt)(t)\ex^{\nu t} - \eps\,(\IO\ft,\Delta\wt)(0)\\
	&+ \eps \int_0^t \bigl[ \nu(\IO\ft,\Delta\wt) + (\dy_\tau\IO\ft,\Delta\wt)
		+ (\IO\ft,\Delta\dstau\wt)\bigr] \ex^{\nu\tau} \dtau,
\end{aligned}\end{equation}
which can be bounded as in the $L^2$ case, giving
\begin{equation}\begin{aligned}
    -&2\int_0^t (\ft,\Delta\wt)\,\ex^{\nu\tau}\dtau
        \le \eps c\,\bigl[|\Delta\ft|\,|\Delta\wt|\bigr](t)\,\ex^{\nu t}
	+ \eps c\,\bigl[|\Delta\ft|\,|\Delta\wt|\bigr](0)\\
	&+ \eps c \int_0^t \bigl\{ (1+\mu)\,|\gb^3\omega|^2 + \mu\,|\gb^3\ft|^2 + |\gb\dy_\tau\ft|^2 + |\gb^3\omega|\,|\gb\omega|\,|\Delta\ft| \bigr\}\,\ex^{\nu\tau}\dtau.
\end{aligned}\end{equation}

For the nonlinear term, we use the fact that $B(\wb,\wb)=0$ to write
\begin{equation}
   (B(\omega,\omega),A\wt) = (B(\wb,\wt),A\wt) + (B(\wt,\wb),A\wt) + (B(\wt,\wt),A\wt),
\end{equation}
and, using $(B(\omega^\sharp,\wt),A\wt) = (B(\gb\omega^\sharp,\wt),\gb\wt)$,
previously used in \eqref{q:whs}, we bound
\begin{equation}\label{q:bdawt}\begin{aligned}
   \bigl|(B(\wt,\wt),A\wt)\bigr|
	&= \bigl|(B(\gb\wt,\wt),\gb\wt)\bigr|
	\le c\,|\wt|_{L^2}^{}|\gb\wt|_{L^4}^2
	\le c\,|\wt|\,|\gb\wt|\,|\Delta\wt|\\
	&\le \frac\mu4\,|\Delta\wt|^2 + \frac{c}{\mu}\,|\wt|^2|\gb\wt|^2\\
   \bigl|(B(\wb,\wt),A\wt)\bigr|
	&= \bigl|(B(\wb',\wt),\dy_y\wt)\bigr|
	\le c\,|\wb|_{L^\infty}^{}|\gb\wt|_{L^2}^2\\
   \bigl|(B(\wt,\wb),A\wt)\bigr|
	&\le |\Delta\wt|_{L^2}^{}|\gb^{-1}\wt|_{L^\infty}^{}|\wb'|_{L^2}^{}
	\le |\Delta\wt|_{L^2}^{}|\gb\wt|_{L^2}^{}|\wb'|_{L^2}^{}\\
	&\le \frac\mu4\,|\Delta\wt|^2
		+ \frac{c}\mu\,|\gb\wt|^2|\wb'|^2.
\end{aligned}\end{equation}
Using Poincar{\'e} inequality on the last term in (\ref{q:bdawt}c), we obtain
\begin{equation}\begin{aligned}
   2\int_0^t (B(\omega,\omega),\Delta\wt)\, &\ex^{\nu\tau}\dtau\\
	&\le c_2(\mu) \int_0^t \Bigl\{\bigl[|\wb|_{L^\infty}^{} + |\gb\omega|^2\bigr]\,|\gb\wt|^2 + \frac\mu2\,|\Delta\wt|^2 \Bigr\}\,\ex^{\nu\tau}\dtau.
\end{aligned}\end{equation}
After moving the $|\Delta\wt|^2$ to the left-hand side,
a factor of $\eps$ can be obtained by pulling the square bracket outside the
integral and using (\ref{q:bdwt}b).
Collecting, we have
\begin{equation}\label{q:bdwt1}\begin{aligned}
   &|\gb\wt(t)|^2 + \frac\mu2 \int_0^t \ex^{\nu(\tau-t)}|\Delta\wt|^2 \dtau
	\le \ex^{-\nu t}|\gb\wt(0)|^2
	+ c\eps\,\sup_{t'>0}\,|\Delta\ft(t')|\,|\Delta\omega(t')|\\
	&\quad {}+ c\eps \int_0^t \bigl\{ (1+\mu)\bigl(|\gb^3\omega|^2 + |\gb^3\ft|^2\bigr) + |\gb\dy_\tau\ft|^2  + |\gb\omega|^2|\Delta\ft|^2 \bigr\}\, \ex^{\nu(\tau-t)}\dtau\\
	&\quad {}+ \eps\,\cnst3(\mu)\,M_0\,\sup_{t'>0}\,\bigl\{|\wb(t')|_{L^\infty}^{} + |\gb\omega(t')|^2\bigr\}.\\
\end{aligned}\end{equation}
Arguing as in the $L^2$ case,
$f\in L_t^\infty H_{\xb}^3$ and $\dy_t\ft\in L_t^\infty H_{\xb}^1$
gives us an $\Oh(\sqrt\eps)$ bound for $\wt(t)$ in $L_t^\infty H_{\xb}^1$
uniform for large $t$.

\medskip
Bounds in $H^s$ can now be obtained inductively.
Assuming that Theorem~\ref{t:os} below holds for $s-1$ (we just showed
that it does for $s=2$), we multiply \eqref{q:abl} by $A^s\wt$
and integrate the resulting equation in time as above.
We bound the nonlinear term
$(B(\omega,\omega),A^s\wt)=(B(\wt,\wt),A^s\wt)+(B(\wb,\wt),A^s\wt)+(B(\wt,\wb),A^s\wt)$
as follows.
The first term is bounded exactly as in \eqref{q:whs}--\eqref{q:interp0},
\begin{equation}\label{q:k01}
   \bigl|(B(\wt,\wt),A^s\wt)\bigr|
	\le \frac\mu4\,|\gb^{s+1}\wt|^2 + \frac{c(s)}\mu\,|\wt|^2|\gb^s\wt|^2.
\end{equation}
We bound the next term by [cf.\ \eqref{q:whs}], with $|\alpha|=s$
and $1\le|\beta|=r\le s$,
\begin{equation}\label{q:k02}\begin{aligned}
   \bigl|(B(\wb,\wt),A^s\wt)\bigr|
	&\le c\,\tssum_{\alpha\beta}\,|\FD^\beta\bar\vb|_{L^\infty}^{}|\FD^{\alpha-\beta}\gb\wt|_{L^2}^{}|\FD^\alpha\wt|_{L^2}^{}\\
	&\le c(s)\,\tssum_{r=1}^s\,|\gb^r\wb|\,|\gb^{s-r+1}\wt|\,|\gb^s\wt|\\
	&\le c(s)\,|\gb^s\wb|\,|\gb^s\wt|^2.
\end{aligned}\end{equation}
Finally, we bound the last term as, where now $0\le|\beta|=r\le s=|\alpha|$,
\begin{equation}\label{q:k03}\begin{aligned}
   \bigl|(B&(\wt,\wb),A^s\wt)\bigr|
	\le c\,\tssum_{\alpha\beta}\,|\FD^\beta\tilde\vb|_{L^4}^{}|\FD^{\alpha-\beta}\gb\wb|_{L^2}^{}|\FD^\alpha\wt|_{L^4}^{}\\
	&\le c(s)\,\tssum_{r=0}^s\,|\wt|_{H^{r-1/2}}^{}\,|\wb|_{H^{s-r+1}}^{}\,|\wt|_{H^{s+1/2}}^{}\\
	&\le \frac\mu4\,|\gb^{s+1}\wt|^2 + c(s,\mu)\,\tssum_{r=0}^s\,|\gb^{r-1}\wt|^{2/3}|\gb^r\wt|^{2/3}|\gb^s\wt|^{2/3}|\gb^{s-r+1}\wb|^{4/3}\\
	&\le \frac\mu4\,|\gb^{s+1}\wt|^2 + c(s,\mu)\,|\gb^s\wt|^2\,|\gb^{s+1}\wb|^{4/3}.
\end{aligned}\end{equation}
Moving the $|\gb^{s+1}\wt|^2$ in \eqref{q:k01} and \eqref{q:k03} to the
left-hand side of the main inequality, the right-hand side depends
at most on $|\gb^s\wt|^2$,
which is of $\Oh(\eps)$ in $L_t^2$ from step $s-1$, and on $|\gb^{s+1}\wb|^2$.
As before, the worst term (i.e.\ that requires the highest derivative on
$f$) in fact comes from bounding $(f,A^s\wt)$.

\vfill\eject

We summarise our results as:

\begin{theorem}\label{t:os}
Let the initial data $\vb(0)\in L^2(\Dom)$ and the forcing be bounded as
\begin{equation}
   K_s(f):= \supt{\gb^{s+2}f} + \supt{\gb^s\dy_tf} < \infty.
\end{equation}
Then there exist $\TB_s(|\vb(0)|_{L^2}^{},K_s;\mu)$ and
$M_s(K_s;\mu)$ such that
\begin{equation}\label{q:bdws}\begin{aligned}
   &|\gb^s\wt(t)|_{L^2}^2 \le \eps\,M_s,\\
   &\mu\,\ex^{-\nu (t+t')}\int_t^{t+t'} \ex^{\nu\tau} |\gb^s\wt(\tau)|^2 \dtau
	\le \eps\,M_s
\end{aligned}\end{equation}
for all $t\ge \TB_s$.
\end{theorem}


\subsection{Higher-order Bounds}\label{s:ho}

As in \cite{temam-dw:lbal}, one can obtain bounds that scale as $\eps^{n/2}$
for $\wt$ when the force $f$ is independent of time;
see \cite{mah:th}.


\section{Stability and the Global Attractor}\label{s:attr}

When the forcing $f$ is independent of time, the existence of the global attractor
$\Attr$ follows, just as for the non-rotating 2d Navier--Stokes equations,
from the uniform long-time bounds in Section~\ref{s:pre},
where the planetary rotation does not appear at all.
In the non-rotating case, the Hausdorff dimension of $\Attr$ is bounded by
\begin{equation}\label{q:dim0}
   \dima \le c\,G^{2/3}(1 + \log G)^{1/3}
\end{equation}
where in our notation the Grashof number is
\begin{equation}
   G := |\gb^{-1} f|_{L^2}^{}/\mu^2.
\end{equation}
The rotation not posing any extra essential difficulty,
the usual analysis, e.g.\ \cite[\S9.2]{doering-gibbon:aanse},
carries over essentially line-by-line to our case,
giving the bound \eqref{q:dim0} also for the rotating case \eqref{q:dwdt}.

As discussed in the introduction, and following our results that the
flow becomes more zonal (``ordered'') as $\eps\to0$, we expect the
dimension of the attractor to decrease as $\eps\to0$.
In this section, we use a simple computation similar to that used
for Theorem~\ref{t:o1} to show that $\dima=0$ for $\eps$
sufficiently small.

\begin{theorem}\label{t:attr}
Let the forcing $f$ be time independent, $\dy_tf=0$, and assume the
hypotheses of Theorem~\ref{t:o1}, i.e.\ $\vb(0)\in L^2(\Dom)$ and
\begin{equation}
   |\gb^2f|_{L^2}^{} < \infty.
\end{equation}
Then there exists an $\eps_*(|\gb^2f|;\mu)$ such that,
for all $\eps<\eps_*$,
\begin{equation}\label{q:dimaz}
   \dima = 0.
\end{equation}
\end{theorem}

Since $\Attr$ is connected, \eqref{q:dimaz} implies that $\Attr$
consists of a single point, that is, a steady flow $\omega_*$ to
which all bounded solutions converge.
Following Theorems \ref{t:o1} and~\ref{t:os}, this steady flow is nearly,
but not exactly, zonal (except in the non-generic case when $\ft=0$).
Heuristically, an approximation to $\omega_*$ is the steady flow
\begin{equation}
   \omega_*^{(1)} = -\mu^{-1}\Delta^{-1}\fb + \eps L^{-1}\ft,
\end{equation}
which satisfies
\begin{equation}
   \frac1\eps L\omega_*^{(1)} + B(\omega_*^{(1)},\omega_*^{(1)}) + \mu A\omega_*^{(1)} = f
\end{equation}
up to $\Oh(\eps)$.
More careful work would be needed to determine $\omega_*$ exactly.

In turbulence parlance, the smallness of $\eps$ demanded by
Theorem~\ref{t:attr} implies that the Rhines scale
\cite{vallis:aofd} is so large that it overwhelms the entire
spectral range, rendering the dynamics trivial.

A general result related to ours is described in
\cite[ch.~18]{chepyzhov-vishik:aemp}, where the trajectory attractor
$\Attr_\epsilon$ of a dynamical system depending on $t/\epsilon$
(formally, in our case $\Attr_\epsilon$ would simply be the attractor
$\Attr$ for $\eps>0$) converges weakly to the attractor $\Attr_0$ of
the corresponding averaged system.
Formally averaging our equations following this construction (which
does not apply directly to our case, in which the oscillations have
an infinite number of frequencies which accumulate at zero),
we obtain purely zonal NSE, whose dynamics is trivial and whose attractor
thus has dimension zero.
This is of course consistent with our results: strong convergence
at finite $\eps$ of $\Attr$ to a point (which becomes zonal as
$\eps\to0$).

\medskip\noindent{\sc Proof.}
Fix a solution $\omega(t)$ of \eqref{q:abl} that lives on $\Attr$,
so the bounds \eqref{q:bdws} hold for all $t$.
We consider a nearby solution $\omega(t)+\phi(t)$.
The linearised evolution equation for $\phi$ is then
\begin{equation}\label{q:dftdt}\begin{aligned}
   \dy_t\phi &= -(\sgb\ilapl\omega)\cdot\gb\phi - (\sgb\ilapl\phi)\cdot\gb\omega(t)
	- \frac1\eps \dy_x\ilapl\phi + \mu\Delta\phi\\
	&= -B(\omega,\phi) - B(\phi,\omega) - \frac1\eps L\phi - \mu A\phi
	=: \Lw(t)\phi.
\end{aligned}\end{equation}
Multiplying this by $\phi$ in $L^2$ and noting that $(B(\omega,\phi),\phi)=0$,
we obtain
\begin{equation}\label{q:dfdt}\begin{aligned}
   \frac12\ddt{\;}|\phi|^2 + \mu\,|\gb\phi|^2 &= (B(\phi,\phi),\omega)\\
	&= (B(\phi,\phi),\wb) + (B(\phi,\phi),\wt).
\end{aligned}\end{equation}

For the first term, we split $\phi=\phb+\pht$ in analogy with $\omega=\wb+\wt$
to get
\begin{equation}
   (B(\phi,\phi),\wb) = (B(\pht,\pht),\wb)
\end{equation}
using the (now familiar) facts that $B(\phb,\phb)=0$ and all
tilde-bar-bar terms vanish.
Using Poincar{\'e} inequality in \eqref{q:dfdt} gives us [cf.\ \eqref{q:dwtdt}]
\begin{equation}
   \ddt{\;}\bigl(\ex^{\nu t}|\phi|^2\bigr) + \mu\ex^{\nu t}|\gb\phi|^2
	\le 2\,\ex^{\nu t}(B(\pht,\pht),\wb) + 2\,\ex^{\nu t}(B(\phi,\phi),\wt),
\end{equation}
which integrates to
\begin{equation}\label{q:l01}\begin{aligned}
   |\phi(t)|^2\,\ex^{\nu t} &+ \mu \int_0^t |\gb\phi|^2\,\ex^{\nu\tau} \dtau\\
	&\le |\phi(0)|^2
	+ 2 \int_0^t \bigl\{ (B(\pht,\pht),\wb) + (B(\phi,\phi),\wt) \bigr\} \,\ex^{\nu\tau} \dtau.
\end{aligned}\end{equation}
We bound the last term of the integrand using
\begin{equation}\begin{aligned}
   (B(\phi,\phi),\wt)
	&\le c\,|\gb^{-1}\phi|_{L^\infty}^{}|\gb\phi|_{L^2}^{}|\wt|_{L^2}^{}\\
	&\le \cnst4\,|\gb\phi|^2|\wt|_{L^2}.
\end{aligned}\end{equation}
The other term needs to be integrated by parts,
\begin{equation}\label{q:l02}\begin{aligned}
   &\int_0^t (B(\pht,\pht),\wb)\,\ex^{\nu\tau}\dtau
	= \eps\,(\BO(\pht,\pht),\wb)(t)\,\ex^{\nu t}
	- \eps\,(\BO(\pht,\pht),\wb)(0)\\
	&\hbox to15pt{} {}- \eps \int_0^t \bigl\{ \nu\,(\BO(\pht,\pht),\wb)
		+ (\BO(\pht,\pht),\dy_\tau\wb)
		+ 2\,(\BO(\dstau\pht,\pht),\wb)
		\bigr\}\,\ex^{\nu\tau}\dtau
\end{aligned}\end{equation}
where $\dy_t^*\phi=-B(\omega,\phi)-B(\phi,\omega)-\mu A\phi$.
We bound the endpoint terms using
\begin{equation}
   2\,\bigl|(\BO(\pht,\pht),\wb)\bigr| \le \cnst5 |\pht|^2|\wb'|_{L^\infty}^{}.
\end{equation}
It remains to bound the integrand in \eqref{q:l02}:
\begin{equation}\begin{aligned}
   \bigl|(\BO(\pht,\pht),\wb)\bigr|
	&\le c\,|\dy_y\pht|_{L^2}^{}|\pht|_{L^4}^{}|\wb|_{L^4}^{}\\
	&\le c\,|\gb\pht|^2|\wb|_{L^4}^{}
\end{aligned}\end{equation}
\begin{equation}\begin{aligned}
   \bigl|(\BO(\pht,\pht),\dy_t\wb)\bigr|
	&\le c\,|\dy_y\pht|_{L^2}^{}|\pht|_{L^{10}}^{}|\dy_t\wb|_{L^{5/2}}^{}\\
	&\le c\,|\gb\pht|^2|\dy_t\wb|_{L^{5/2}}^{}
\end{aligned}\end{equation}
Recalling \eqref{q:dftdt} for the last term in \eqref{q:l02}, we bound
\begin{align*}
   \bigl|(\BO(\tilde B(\phi,\omega),\pht),\wb)\bigr|
	&\le c\,|\gb^{-1}\phi|_{L^\infty}^{}|\gb\omega|_{L^2}^{}|\pht|_{L^{10}}^{}|\wb'|_{L^{5/2}}^{}\\
	&\le c\,|\gb\phi|^2|\wb'|_{L^{5/2}}^{}|\gb\omega|_{L^2}^{}
\end{align*}
\begin{align}
   \bigl|(\BO(\tilde B(\omega,\phi),\pht),\wb)\bigr|
	&\le c\,|\gb^{-1}\omega|_{L^\infty}^{}|\gb\phi|_{L^2}^{}|\pht|_{L^{10}}^{}|\wb'|_{L^{5/2}} \notag\\
	&\le c\,|\gb\phi|^2|\wb'|_{L^{5/2}}^{}|\gb\omega|_{L^2}^{}
\end{align}
\begin{align*}
   \bigl|(\BO(\Delta\pht,\pht),\wb)\bigr|
	&= \bigl|(\BO(\gb\pht,\gb\pht),\wb)\bigr| + \bigl|(\BO(\dy_y\pht,\pht),\wb')\bigr|\\
	&\le c\,|\gb\pht|^2|\wb'|_{L^\infty}^{} + c\,|\gb\pht|_{L^2}^{}|\pht|_{L^{10}}^{}|\wb''|_{L^{5/2}}^{}\\
	&\le c\,|\gb\pht|^2|\wb''|_{L^{5/2}}^{}.
\end{align*}
Collecting, \eqref{q:l01} now implies
\begin{equation}\label{q:a0}\begin{aligned}
   |\phi(t)|^2\bigl(1 - \eps\cnst5|\wb'(t)|_{L^\infty}^{}\bigr)
	+ \int_0^t \bigl\{ \mu - \eps\,N(\tau) &- \cnst4|\wt(\tau)|_{L^2}^{} \bigr\}\,|\gb\pht|^2\,\ex^{\nu(\tau-t)}\dtau\\
	&\le \ex^{-\nu t}|\phi(0)|^2\bigl(1 + \eps\cnst5|\wb'(0)|_{L^\infty}^{}\bigr)
\end{aligned}\end{equation}
where
\begin{equation}
   N(t) := \cnst6 \bigl\{ \mu\,|\wb''|_{L^{5/2}}^{}
		+ |\wb'|_{L^{5/2}}^{}|\gb\omega|_{L^2}^{}
		+ |\dy_t\wb|_{L^{5/2}}^{} + |\wb|_{L^4}^{} \bigr\}(t).
\end{equation}
By \eqref{q:bdwm}, $f\in H^2$ implies that $\omega\in H^3$ with
uniform bound in $t$ since we are already on the attractor,
and by Theorem~\ref{t:o1} we can find an $\eps_*$ so small that,
for $\eps<\eps_*$,
\begin{equation}
   \sup_{t>0}\, \bigl\{ \eps N(t) + \cnst4|\wt(t)|_{L^2}^{} \bigr\} < \mu.
\end{equation}
Requiring furthermore that $\eps_*$ also implies, for all $\eps<\eps_*$,
\begin{equation}
  \eps\cnst5\,\sup_{t>0}\,|\wb'(t)|_{L^\infty}^{} < 1.
\end{equation}
These and \eqref{q:a0} then imply that
\begin{equation}
   |\phi(t)|^2 \le C(\cdots)\,\ex^{-\nu t}|\phi(0)|^2,
\end{equation}
in other words, all phase space volumes contract and thus
the global attractor has dimension zero.

It is clear from the above proof that our solution $\omega(t)$
is linearly stable.
Since \eqref{q:dftdt} only differs by $B(\phi,\phi)$ from the
nonlinear system,
the fact that $(B(\phi,\phi),\phi)=0$ implies that stability
also holds under the same hypotheses for the full nonlinear system.



\appendix

\section{An $L^\infty$ Inequality}\label{s:linfty}

\begin{lemma}
Let $u$ and $v\in H^2(\Dom)$ have zero integrals and are $L^2$ orthogonal,
\begin{equation}
   (u,v)_{L^2}^{} = 0,
\end{equation}
and let $w=u+v$.
Then the following Agmon inequality holds,
\begin{equation}\label{q:agmon}
   |u|_{L^\infty}^{} \le c\,|\gb w|\,\Bigl(\log\frac{|\Delta w|}{\cpoi|\gb w|} + 1\Bigr)^{1/2}.
\end{equation}
\end{lemma}

Before the proof, we note that the interpolation inequality
\begin{equation}
   |\gb w|^2 \le \cnst9\,|w|\,|\Delta w|
\end{equation}
can be written as
\begin{equation}\begin{aligned}
   &2\log|\gb w| \le \log|w| + \log|\Delta w| + \log\cnst9\\
   \Leftrightarrow\qquad
   &\log|\gb w| - \log|w| \le \log|\Delta w| - \log|\gb w| + \log\cnst9\\
   \Leftrightarrow\qquad
   &\log\frac{|\gb w|}{\cpoi|w|} \le \log\frac{|\Delta w|}{\cpoi|\gb w|} + \log\cnst9,
\end{aligned}\end{equation}
which can be used to simplify, e.g., $|w|_{L^\infty}^{}|\gb w|_{L^\infty}^{}$
when bounded using \eqref{q:agmon}.

\medskip\noindent{\sc Proof.}
For most of this proof, up to \eqref{q:kappadef} below, we follow
\cite[Lemma~7.1]{doering-gibbon:aanse} exactly.
For conciseness, we put $L_1=L_2=1$ but keep the Poincar{\'e} constant $\cpoi$.
With $\kappa>0$, we expand $u$ in Fourier series
\begin{equation}
   u(\xb) = \tssum_{|\kb|<\kappa}\,u_\kb\ex^{\im\kb\cdot\xb}
	+ \tssum_{|\kb|\ge\kappa}\,u_\kb\ex^{\im\kb\cdot\xb}
	=: u^<(\xb) + u^>(\xb),
\end{equation}
and analogously for $v$ and $w$.
Then
\begin{equation}\begin{aligned}
   |u|_{L^\infty}^{}
	&= \textstyle\sup_{\xb}\,\bigl|\tssum_\kb\, u_\kb\ex^{\im\kb\cdot\xb}\bigr|
	\le \tssum_{|\kb|<\kappa}\,|u_\kb^{}| + \tssum_{|\kb|\ge\kappa}\,|u_\kb^{}|\\
	&=: \tssum^<\,|\kb|^{-1}|\kb|\,|u_\kb^{}|
		+ \tssum^>\,|\kb|^{-2}|\kb|^2|u_\kb^{}|\\
	&\le \Bigl(\tssum^<\,|\kb|^{-2}\Bigr)^{1/2}\Bigl(\tssum^<\,|\kb|^2|u_\kb^{}|^2\Bigr)^{1/2}\\
	&\hbox to 50pt{} + \Bigl(\tssum^>\,|\kb|^{-4}\Bigr)^{1/2}\Bigl(\tssum^>\,|\kb|^4|u_\kb^{}|^2\Bigr)^{1/2}
\end{aligned}\end{equation}
Now on the right-hand side,
$\sum^<|\kb|^{-2}\le c\,\log\kappa$ and $\sum^>|\kb|^{-4}\le c/\kappa^2$,
so fixing
\begin{equation}\label{q:kappadef}
   \kappa = |\Delta w|/(\cpoi|\gb w|),
\end{equation}
the lemma follows from
\begin{equation}\begin{aligned}
   |u|_{L^\infty}^{}
	&\le c\,|\gb u^<|\Bigl(\log\frac{|\Delta w|}{\cpoi|\gb w|}\Bigr)^{1/2}
	+ c\,|\Delta u^>|\,\frac{|\gb w|}{|\Delta w|}\\
	&\le c\,|\gb w|\Bigl(\log\frac{|\Delta w|}{\cpoi|\gb w|}\Bigr)^{1/2}
	+ c\,|\gb w|.
\end{aligned}\end{equation}




\end{document}